\documentclass[12pt,british]{elsarticle}
\usepackage[T1]{fontenc}
\usepackage[latin1]{inputenc}
\usepackage{textcomp}
\usepackage{mathrsfs}
\usepackage{amsmath}
\usepackage{amssymb}
\usepackage{graphicx}
\usepackage{setspace}
\makeatletter

\newcommand{\noun}[1]{\textsc{#1}}
\providecommand{\tabularnewline}{\\}

\makeatother

\usepackage{babel}
\begin{document}

\title{A general stability criterion for switched linear systems having
stable and unstable subsystems}

\author{Jesús San Martín}

\address{EUITI Universidad Politécnica de Madrid, Ronda de Valencia, 3. 28012
Madrid (Spain)\\
Dep. Física Matemática y de Fluidos, UNED, Pº Senda del Rey, 9. 28040
Madrid (Spain)}

\ead{jsm@dfmf.uned.es}

\author{Anthony G. O'Farrell}

\address{Dept. of Mathematics and Statistics, NUI Maynooth. Co. Kildare, Ireland}

\ead{admin@maths.nuim.ie}
\begin{abstract}
We report conditions on a switching signal that
guarantee that solutions of a switched linear systems converge asymptotically
to zero. These conditions are apply to continuous, discrete-time
and hybrid switched linear systems, both those having stable subsystems
and mixtures of stable and unstable subsystems. \end{abstract}
\maketitle

\section{Introduction}

In Science and Engineering, systems are frequently met that consist
of a family of subsystems and a switching signal which determines
which subsystem is activated each time.

When all the subsystems are linear, one has a switched linear system
\begin{equation}
\overset{.}{x}(t)=A_{\sigma(t)}x(t)\label{eq:(1)}
\end{equation}
where $\sigma:\left[0,+\infty\right)\rightarrow\left\{ 1,\cdots,n\right\} $
is the switching signal and $A_{i}:\mathbb{R}^{m}\rightarrow\mathbb{R}^{m}\quad(i=1,\cdots,n)$
are matrices that characterise the subsystems. We will assume the
condition that there are a finite number of switches in each finite
time interval \citep{Liberzon1999,Branicky1998}.

The great number of areas in which switched linear systems appear
makes their study a matter of real concern and great importance \citep{Liberzon1999}.
Its theoretical importance \citep{Liberzon1999,Sun2005,Sun2005a}
derives from its practical importance (power systems, control process,
automotive industry, mechanical systems): one needs to understand
under what circumstances the system (\ref{eq:(1)}) is stable, or
what switching signals make the systems stable (see the survey \citep{Lin2009}
for more details). Many times, time-delay must be taken into account
\citep{Wang2011,Qiu2009} because it plays an important role in many
practical systems, such as chemical processes, nuclear reactors, engines,
and so on \citep{Guo2010,Orlov2009}.

Liberzon and Morse \citep{Liberzon1999} formulated three basic problems
in relation to the stability of switched systems.

\textsl{{}``Problem A}\textit{: Find conditions that guarantee that
the switched system is asymptotically stable for any switching signal''.}

\textsl{{}``Problem B}\textit{: Identify those classes of switching
signals for which the switched system is asymptotically stable''.}

\textit{{}``Problem C: Construct a switching signal that makes the
switched system asymptotically stable''.}

The condition of asymptotic stability referred to \textit{Problem
A}, is desirable in practical applications. The theorems that provide
solution (or partial solutions) to in \textit{Problem A} involve conditions
that can be restrictive for applications: the existence of Lyapunov
functions, symmetric systems, pairwise commutativity of the subsystems,
Lie-algebraic conditions, \ldots{} \citep{Branicky1998,Sun2005,Liberzon1999a,Benzaouia2011}.
On the other hand, it is well-known that there exist systems that
exhibit instability even though all their subsystems are asymptotically
stable \citep{Liberzon1999,Branicky1998}. As a result, one sees the
necessity of solving \textsl{Problem B} in practice, in order to deal
with the applications. Usually, \textsl{Problem B} is studied under
the assumption that all the individual subsystems are asymptotically
stable \citep{Liberzon1999,Liberzon1999a}. However, for some applications
it is convenient to allow subsystems that may be both stable and unstable
(for instance, unstable subsystems have to be considered when a failure occurs).

In this paper, we establish conditions on the switching signal of
a switched linear system that are sufficient to ensure asymptotic
stability. We allow both stable and unstable subsystems in the switched
linear system \citep{Zhai2000}. Furthermore this switched linear
system has no restrictive conditions involved, such as the existence
of Lyapunov functions, symmetries, pairwise commutativity of the subsystems,
Lie-algebraic conditions, etc.

In the following lines we are putting our results in context and indicate
which are the novel contributions. As we said above, switched systems
with all the subsystems stable can become unstable for certain switching
signals. The contrary is also true, that switched systems with some
of the subsystems unstable can become stable for certain switching
signals \citep{Decarlo2000,Liberzon1999a}. From what we have just
said it is plain to see that the stability of switched systems depends
on both the dynamics of each subsystem and the switching signals.

Intuitively it is easy to imagine a switching signal that maintains
stable a switched linear system with all its subsystems stable. It
is enough to stay for a long time in one of the subsystems and jump
into another with a low frequency. This idea, and the mathematical
concept of \textquotedblleft{}low frequency of jumping between subsystems\textquotedblright{},
was developed by Morse and Hespanha with the concepts of dwell time
and average dwell time switching \citep{Hespanha2004,Hespanha1999}.

These results about stability were extended in \citep{Zhai2000} to
take into account unstable subsystems, that is switched systems having
 both stable and unstable
subsystems. The underlying idea is very
similar to the one just described above. Stay for long times in stable
subsystems in such a way that the total activation time of stable
subsystems is relatively big compared with that of unstable ones,
 and furthermore
the system should spend long times in every stable subsystems (low frequency
of jumping).

Dwell time and average dwell time switching have proved to be very
flexible and powerful tools to determine what conditions switching
signal should satisfy to guarantee the stability of switched systems.
As a result, there is a fruitful line of research since the pioneering
ideas to the present day \citep{Zhang2009_a,Zhao2012,Lian2011,Zhang2009_b}.
New approaches using fuzzy time control have also been taken \citep{Wu2011}.

The results shown in this paper also require that stable subsystems
control the unstable one, but in a different way. Average dwell time
switching requires that the average interval between consecutive switchings
must be bigger than a constant, although consecutive switchings can
be separated by less than that constant. In contrast to this idea
the theorems we show do not need time intervals, we use the number
of times that the system switches. We are speaking of a discrete variable
in contrast to a continuous one. That is not the only difference.
We require that stable subsystems control unstable subsystems 
pair by pair,
as we will explain below. A consequence is that a switched system
can spend more time in unstable subsystems than in stable ones and
even then the switched system be stable, as is illustrated in example
2.

The paper is organised as follows. First we work on continuous switched
linear systems, then on discrete systems. Afterwards we combine these
to study hybrid systems. Then we apply our results about problems
of \textit{type} \textit{B} to the design of switching signals in
order to solve problems of \textit{type C}.

\section{Definitions}

The following definitions are necessary for further discussion.

$t_{i_{j}}$ denotes the time during which the system (\ref{eq:(1)})
is ruled by matrix $A_{i}$, when this matrix is switched to for the
$j$-th time.

$m_{i}(t):[0,\infty)\rightarrow\mathbb{N}$ denotes how many times
the matrix $A_{i}$ has been switched up to time $t$ .

\section{Continuous systems}

We first investigate continuous systems, then we will reformulate
the results of this section to approach discrete systems. 

The next theorem gives us sufficient conditions to get asymptotic
stability in a continuous system having both stable and unstable subsystems.
It is important to emphasise two things. First, the system will be
asymptotically stable even though it has unstable subsystems. Second,
we are not using average dwell time switching. This theorem will allow
us the stabilisation of a switched linear system by a switching strategy.

\paragraph{Theorem~1 \label{thm:1}}

Consider a switched linear system of the form (\ref{eq:(1)}) that
satisfies:
\begin{enumerate}
\item [(i)] $m_{i}(t)\rightarrow\infty\quad\mbox{as}\quad t\rightarrow\infty\quad\mbox{for}\quad i=1,\cdots,n$
\item [(ii)]$<e^{A_{i}}>=\left[\prod_{i_{j}=1}^{m_{i}}\left\Vert e^{A_{i}t_{i_{j}}}\right\Vert \right]^{\frac{1}{m_{i}}}\rightarrow c_{i}<+\infty\quad\mbox{as}\quad m_{i}\rightarrow+\infty$
where

\[
\begin{array}{cc}
c_{i}\leq1-\varepsilon_{\beta}, & \quad\varepsilon_{\beta}>0,\qquad i=1,\cdots,k-1\\
c_{i}\geq1, & i=k,\cdots,n
\end{array}
\]

\item [(iii)] $\frac{m_{j}(t)}{m_{\rho}(t)}\rightarrow k_{j,\rho}\geq1\quad\mbox{as}\quad t\rightarrow+\infty\quad\mbox{for}\quad\begin{array}{c}
j=1,\cdots,k-1\\
\rho=k,\cdots,n
\end{array}$
\item [(iv)] $\prod_{i=1}^{n}c_{i}\leqslant1-\varepsilon_{\alpha}\quad\varepsilon_{\alpha}>0$
\end{enumerate}
Then each solution $x(t)$ of (\ref{eq:(1)}) tends asymptotically
to zero.

\paragraph{Significance}
\begin{enumerate}
\item [i)] Each matrix $A_{i}$ is used infinitely many times. The process
may be random or deterministic. The systems does not stay in any one
state as $t\rightarrow\infty$.
\item [ii)] The geometric average of the norms of the flows for each matrix
$A_{i}$ of (\ref{eq:(1)}) is finite. The matrices $A_{j},\quad j=1,\cdots,k-1$
contract the flow, whereas the matrices $A_{\rho},\quad\rho=k,\cdots,n$
expand it.
\item [iii)] Each contracting matrix dominates each expanding matrix.
\item [iv)] The contracting flows dominate the expanding ones.
\end{enumerate}

\paragraph{Proof}

After a time $t$, the transfer matrix $A_{i}\quad i=1,\cdots,n$
will have been used $m_{i}(t)$ times. Taking norms it results

\begin{equation}
\left\Vert x(t)\right\Vert \leq\prod_{i_{n}=1}^{m_{n}}\left\Vert e^{A_{n}t_{i_{n}}}\right\Vert ...\prod_{i_{1}=1}^{m_{1}}\left\Vert e^{A_{1}t_{i_{1}}}\right\Vert \left\Vert x_{0}\right\Vert \label{eq:(2)}
\end{equation}
 it follows from {\footnotesize (ii) }that
\[
\forall\varepsilon_{i}>0\quad\exists\; N_{\varepsilon_{i}}/m_{i}>N_{\varepsilon_{i}}\rightarrow\left|<e^{A_{i}}>-c_{i}\right|\leqslant\varepsilon_{i}\quad i=1,\cdots,n
\]

Let $N_{\varepsilon}=\underset{i=1,\cdots,n}{\max}\left\{ N_{\varepsilon_{i}}\right\} \quad\varepsilon=\underset{i=1,\cdots,n}{\max}\left\{ \varepsilon_{i}\right\} $.

If

\begin{equation}
m_{i}>N_{\varepsilon}\quad i=1,\cdots,n\label{eq:(3)}
\end{equation}
it follows from (\ref{eq:(2)}) that 
\[
\left\Vert x(t)\right\Vert \leq\left(c_{n}+\varepsilon\right)^{m_{n}}...\left(c_{k}+\varepsilon\right)^{m_{k}}\left(c_{k-1}+\varepsilon\right)^{m_{k-1}}...\left(c_{1}+\varepsilon\right)^{m_{1}}\leq
\]
\begin{equation}
\leq\left(c_{n}+\varepsilon\right)^{M}...\left(c_{k}+\varepsilon\right)^{M}\left(c_{k-1}+\varepsilon\right)^{m}...\left(c_{1}+\varepsilon\right)^{m}\left\Vert x_{0}\right\Vert \label{eq:(4)}
\end{equation}
where $M=\underset{l=k,\cdots,n}{\max}\left\{ m_{l}\right\} ,\quad m=\underset{j=1,\cdots,k-1}{\min}\left\{ m_{j}\right\} $
and $\varepsilon_{i}$ has been chosen so that $\varepsilon_{i}<\varepsilon_{\beta}\quad i=1,...,n$ 

It follows from {\footnotesize (iii)} that
\[
\exists t_{L}/t>t_{L}\rightarrow m\geq M
\]
if 
\begin{equation}
t>t_{L}\label{eq:cinco}
\end{equation}
 it follows from (\ref{eq:(4)}) that
\begin{equation}
\begin{array}{c}
\left\Vert x(t)\right\Vert \leq\left[\left(c_{n}+\varepsilon\right)...\left(c_{1}+\varepsilon\right)\right]^{m}\left\Vert x_{0}\right\Vert =\\
=\left(c_{1}...c_{n}+k\varepsilon\right)^{m}\left\Vert x_{0}\right\Vert \leq\left(1-\varepsilon_{\alpha}+k\varepsilon\right)^{m}\left\Vert x_{0}\right\Vert 
\end{array}\label{eq:(6)}
\end{equation}
where {\footnotesize (iv)} has been used in the last equality.

We choose $\varepsilon<\min\left\{ \varepsilon_{\beta},\frac{\varepsilon_{\alpha}}{k}\right\} $
then $1-\varepsilon_{\alpha}+k\varepsilon<1$

It follows from {\footnotesize (i)} that
\[
\forall N\quad\exists t_{N}/t>t_{N}\rightarrow m_{i}\geq N\quad i=1,\cdots,n
\]
it is sufficient to take $t>\max\left\{ t_{L},\, t_{N_{\varepsilon}}\right\} $
so that (\ref{eq:(3)}) and (\ref{eq:cinco}) are met and consequently
(\ref{eq:(6)}) with $1-\varepsilon_{\alpha}+k\varepsilon<1$. It
follows from {\footnotesize (i)} that $m\underset{t\rightarrow\infty}{\rightarrow\infty}$.
Therefore

\[
\left\Vert x(t)\right\Vert \leq\left(1-\varepsilon_{\alpha}+k\varepsilon\right)^{m}\left\Vert x_{0}\right\Vert \underset{t\rightarrow\infty}{\rightarrow}0
\]
$\blacksquare$

Presented below is a numerical example in the continuous-time domain.

\paragraph{Example 1}

Consider the switched linear system 

\begin{equation}
\overset{.}{x}(t)=A_{\sigma(t)}x(t)\label{eq:ex_1}
\end{equation}

\[
\sigma:\left[0,+\infty\right)\rightarrow\left\{ 1,\cdots,4\right\} 
\]

where switching sequence is generated by random numbers and

~

$A_{1}=\left(\begin{array}{cc}
-2 & 3\\
-4 & 5
\end{array}\right)$ 

~

$A_{2}=\left(\begin{array}{rr}
0 & 6\\
-1 & 5
\end{array}\right)$ 

~

$A_{3}=\left(\begin{array}{rr}
-\frac{29}{4} & \frac{4}{9}\\
-\frac{10}{9} & -\frac{7}{9}
\end{array}\right)$ 

~

$A_{4}=\left(\begin{array}{rr}
-\frac{69}{17} & \frac{2}{17}\\
-\frac{9}{17} & -\frac{50}{17}
\end{array}\right)$ 

~

The subsystems ruled by $A_{1}$ and $A_{2}$ are unstable because
their eigenvalues are respectively $\lambda(A_{1})=\left\{ 1,\,2\right\} $
and $\lambda(A_{2})=\left\{ 2,\,3\right\} $. Whereas $A_{3}$ and
$A_{4}$ , with eigenvalues $\lambda(A_{3})=\left\{ -1,\,-3\right\} $
and $\lambda(A_{4})=\left\{ -3,\,-4\right\} $ respectively, determine
stable subsystems. 

In order to calculate the evolution of the system (\ref{eq:ex_1})
the following expressions are needed.

~

$e^{A_{1}t}=\left(\begin{array}{cc}
-1 & -\frac{3}{4}\\
-1 & -1
\end{array}\right)\left(\begin{array}{cc}
e^{t} & 0\\
0 & e^{2t}
\end{array}\right)\left(\begin{array}{rr}
-4 & 3\\
4 & -4
\end{array}\right)$

~

$e^{A_{2}t}=\left(\begin{array}{cc}
3 & 2\\
1 & 1
\end{array}\right)\left(\begin{array}{cc}
e^{2t} & 0\\
0 & e^{3t}
\end{array}\right)\left(\begin{array}{rr}
1 & -2\\
-1 & 3
\end{array}\right)$

~

$e^{A_{3}t}=\left(\begin{array}{cc}
1 & 2\\
5 & 1
\end{array}\right)\left(\begin{array}{cc}
e^{-t} & 0\\
0 & e^{-3t}
\end{array}\right)\left(\begin{array}{rr}
-\frac{1}{9} & \frac{2}{9}\\
\frac{5}{9} & -\frac{1}{9}
\end{array}\right)$

~

$e^{A_{4}t}=\left(\begin{array}{cc}
1 & 2\\
9 & 1
\end{array}\right)\left(\begin{array}{cc}
e^{-3t} & 0\\
0 & e^{-4t}
\end{array}\right)\left(\begin{array}{rr}
-\frac{1}{17} & \frac{2}{17}\\
\frac{9}{17} & -\frac{1}{17}
\end{array}\right)$

~

After 100000 switches, the result is 

~

$\left\langle e^{A_{1}}\right\rangle =\left[\prod_{1_{j}}^{m_{1}}\left\Vert e^{A_{1}t_{_{1_{j}}}}\right\Vert _{2}\right]^{\frac{1}{m_{1}}}\underset{m_{1}\rightarrow\infty}{\rightarrow}c_{1}=1.106\geq1$

~

$\left\langle e^{A_{2}}\right\rangle =\left[\prod_{2_{j}}^{m_{2}}\left\Vert e^{A_{2}t_{_{2_{j}}}}\right\Vert _{2}\right]^{\frac{1}{m_{2}}}\underset{m_{2}\rightarrow\infty}{\rightarrow}c_{2}=1.198\geq1$

~

$\left\langle e^{A_{3}}\right\rangle =\left[\prod_{3_{j}}^{m_{3}}\left\Vert e^{A_{3}t_{3_{j}}}\right\Vert _{2}\right]^{\frac{1}{m_{3}}}\underset{m_{3}\rightarrow\infty}{\rightarrow}c_{3}=0.985\leq1-\varepsilon_{\beta}\quad\varepsilon_{\beta}>0$

~

$\left\langle e^{A_{4}}\right\rangle =\left[\prod_{4_{j}}^{m_{4}}\left\Vert e^{A_{4}t_{_{4j}}}\right\Vert _{2}\right]^{\frac{1}{m_{4}}}\underset{m_{4}\rightarrow\infty}{\rightarrow}c_{4}=0.101\leq1-\varepsilon_{\beta}\quad\varepsilon_{\beta}>0$

~

with

~

$\prod_{i=1}^{4}c_{i}=0.1321-\varepsilon_{\alpha}\quad\varepsilon_{\alpha}>0$ 

\begin{figure}
\includegraphics[width=0.55\textwidth]{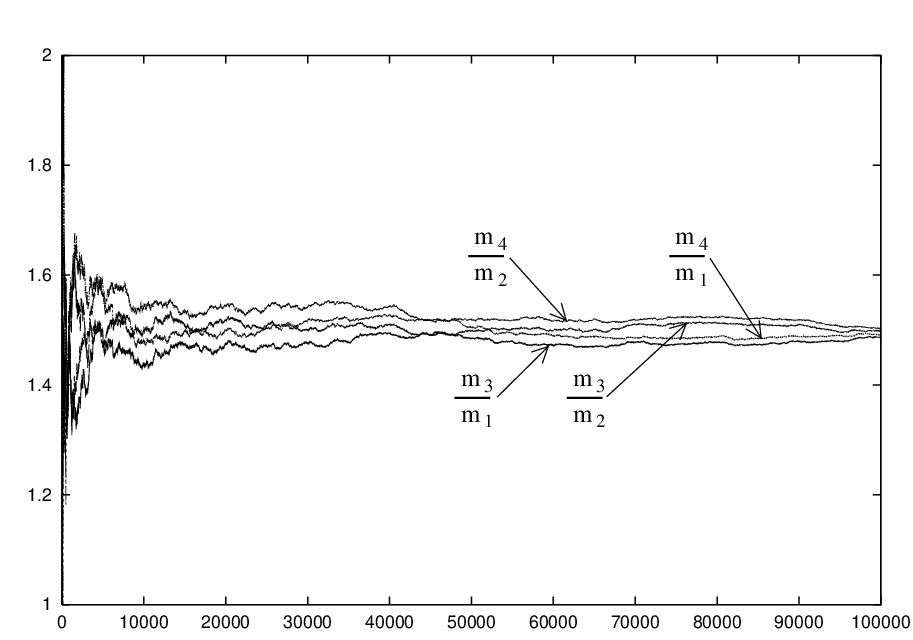}

\caption{\label{fig:1}Temporal evolution of $m_{3}/m_{1}$ , $m_{3}/m_{2}$,
$m_{4}/m_{1}$ and $m_{4}/m_{2}$.}
\end{figure}

~

and (see Figure \ref{fig:1})

~

$\frac{m_{3}}{m_{1}}\underset{t\rightarrow\infty}{\rightarrow}k_{3,1}=1.486$

~

$\frac{m_{3}}{m_{2}}\underset{t\rightarrow\infty}{\rightarrow}k_{3,2}=1.498$

~

$\frac{m_{4}}{m_{1}}\underset{t\rightarrow\infty}{\rightarrow}k_{4,1}=1.492$

~

$\frac{m_{4}}{m_{2}}\underset{t\rightarrow\infty}{\rightarrow}k_{4,2}=1.503$

~

\begin{figure}
\begin{tabular}{cc}
{\footnotesize (a)}\includegraphics[width=0.55\textwidth]{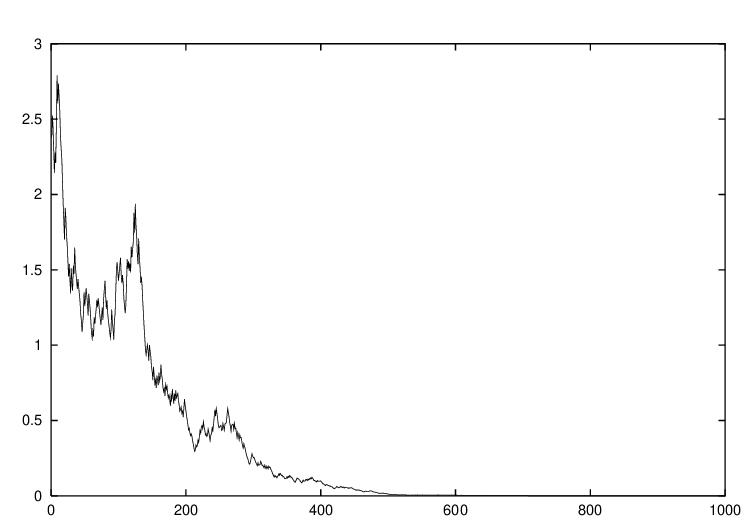} & {\footnotesize (b)}\includegraphics[width=0.55\textwidth]{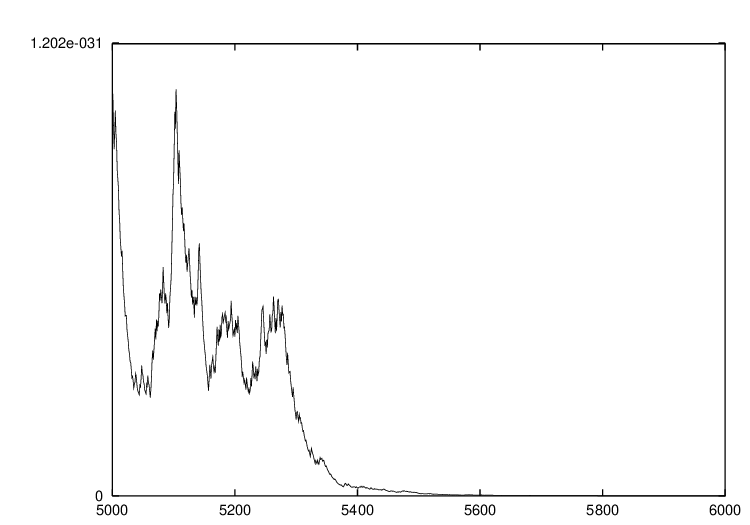}\tabularnewline
\end{tabular}

\caption{\label{fig:2}Temporal evolution of the asymptotic convergence in
norm to zero of continuous system (\ref{eq:ex_1}).}
\end{figure}

So theorem conditions from (i) to (iv) are satisfied and the system
must asymptotically tend to zero, as it is shown in Figures \ref{fig:2}a
and \ref{fig:2}b.

\section{Discrete-time systems}

When time is discrete instead of continuous we have a switched linear
discrete-time system, and the system (\ref{eq:(1)}) is turned into
\begin{equation}
x(n+1)=A_{\sigma(t)}x(n)\label{eq:(b)}
\end{equation}
where $\sigma(t)$ and $A_{\sigma(t)}$ have the same meaning as
in system (\ref{eq:(1)}).

Discrete-time systems are as useful in engineering as continuous-time
systems, and theoretical research is also very active. Furthermore,
they appear in other areas where continuous systems are not found,
as a result of using the transfer matrix method to solve differential
equations \citep{Khorasani2003}. Lately, they are becoming more important
in the study of structures consisting of stiffened plates (naval architecture,
bridge engineering, aircraft design buildings) \citep{Xie2000} and
spatially periodic structures (satellite antennae, satellite solar
panels) \citep{Xie2001}. The theorem, stated some lines below, will
indicate to a designer how to insert panels (given by $A_{i}$ in (\ref{eq:(b)}))
so that oscillations fade off and do not damage the structure.

The theorem 1 can be reformulated for discrete-time systems in the
following way:

\paragraph{Theorem 2 \label{thm:2}}

Consider a switched linear discrete-time of the form (\ref{eq:(b)})
such that
\begin{enumerate}
\item [i)] $m_{i}(t)\underset{t\rightarrow\infty}{\rightarrow}\infty\quad i=1,\cdots,n$
\item [ii)] $\frac{m_{j}(t)}{m_{\rho}(t)}\underset{t\rightarrow\infty}{\rightarrow}k_{j,\rho}\geq1\quad\begin{array}{c}
j=1,\cdots,k-1\\
\rho=k,\cdots,n
\end{array}$
\item [iii)] $\begin{array}{ccc}
\left\Vert A_{j}\right\Vert <1 & \quad & j=1,\cdots,k-1\\
\left\Vert A_{\rho}\right\Vert \geq1 & \quad & \rho=k,\cdots,n
\end{array}$
\item [iv)] $\prod_{i=1}^{n}\left\Vert A_{i}\right\Vert \leq1-\varepsilon_{\alpha}\quad\varepsilon_{\alpha}>0$
\end{enumerate}
Then each solution of \eqref{eq:(b)} tends to $0$ as $n\to+\infty$.

\paragraph{Proof}

After a time $t$ has elapsed, the transfer matrix $A_{i}\quad i=1,\cdots,n$,
will have been used $m_{i}(t)$ times. Hence, taking norms we have
\[
\left\Vert x(t)\right\Vert \leq\prod_{i_{n}=1}^{m_{n}}\left\Vert A_{n}\right\Vert ...\prod_{i_{1}=1}^{m_{1}}\left\Vert A_{1}\right\Vert \left\Vert x_{0}\right\Vert \leq
\]
\begin{equation}
\leq\left(\left\Vert A_{n}\right\Vert ...\left\Vert A_{k}\right\Vert \right)^{M}\left(\left\Vert A_{k-1}\right\Vert ...\left\Vert A_{1}\right\Vert \right)^{m}\left\Vert x_{o}\right\Vert \label{eq:(2-1)}
\end{equation}
where $M=\underset{l=k,\cdots,n}{\max}\left\{ m_{l}\right\} ,\quad m=\underset{j=1,\cdots,k-1}{\min}\left\{ m_{j}\right\} $
and {\footnotesize (iii)} has been used.

It follows from $(ii)$ that $\exists\quad t_{L}/t>t_{L}\rightarrow m\geq M$
.

If 
\begin{equation}
t>t_{L}\;\label{eq:dos-dos}
\end{equation}
it follows from (\ref{eq:(2-1)}) that

\[
\left\Vert x(t)\right\Vert \leq\left(\left\Vert A_{n}\right\Vert ...\left\Vert A_{1}\right\Vert \right)^{m}\left\Vert x_{o}\right\Vert \leq\left(1-\varepsilon_{\alpha}\right)^{m}\left\Vert x_{o}\right\Vert \underset{t\rightarrow\infty}{\rightarrow}0
\]
where {\footnotesize (ii)} has been used and condition (\ref{eq:dos-dos})
is satisfied when $t\rightarrow\infty$. $\blacksquare$

If we consider the problem mentioned at the beginning of the section,
and imagine a {}``solar panel'' with many sections suffering unstable
oscillations then the theorem will indicate the possibility of inserting
a panel to extinguish the vibrations.

Giving that the solar panel is a periodical structure, such that the
switching to its different components would be ruled by a travelling
wave it follows that the switching signal $\sigma(t)$ would be determined
by a deterministic expression; hence, the engineer will have to choose
the materials in the solar panel so that theorem 2 is satisfied and
the travelling wave in it is extinguished. 

A similar argument would allow one to deduce whether a wave would extinguish
in a system ruled by Schrödinger or Maxwell equations \citep{Mayer1999}.

\paragraph{Remark 1}

It is very easy for engineers to decide whether condition $(iv)$ is
satisfied. Given that $\left\Vert A\right\Vert _{2}=\sqrt{\lambda_{max}(A^{*}A)}\quad$
where $\lambda_{max}(A^{*}A)=\max\left\{ \mbox{eigenvalues\,\ of\,}A^{*}A\right\} $
it follows that condition $(iv)$ is satisfied if

\[
\left\Vert A_{n}\right\Vert _{2}...\left\Vert
 A_{1}\right\Vert _{2}=\prod_{i=1}^{n}\sqrt{\lambda_{max}(A_{i}^{*}A)}<1
\]

Below is a numerical example for discrete-time systems.

\paragraph{Example 2}

Consider the switched linear system 

\begin{equation}
x(n+1)=A_{\sigma(t)}x(n)\label{eq:ex_2}
\end{equation}
where switching sequence is generated by random numbers and

~

$A_{1}=\left(\begin{array}{cc}
0 & -\frac{29}{50}\\
\frac{28}{50} & \frac{1}{10}
\end{array}\right)$ 

~

$A_{2}=\left(\begin{array}{cc}
\frac{51}{40} & -\frac{1}{40}\\
-\frac{3}{40} & \frac{49}{40}
\end{array}\right)$ 

~

$A_{3}=\left(\begin{array}{cc}
\frac{57}{50} & \frac{1}{25}\\
\frac{3}{50} & \frac{29}{25}
\end{array}\right)$ 

~

$\left\Vert A_{1}\right\Vert _{2}\simeq0.623<1$

~

with eigenvalues $\lambda(A_{1})=\left\{ \frac{1}{20}+\frac{\sqrt{3223}i}{100},\,\frac{1}{20}-\frac{\sqrt{3223}i}{100}\right\} $ 

~

$\left\Vert A_{2}\right\Vert _{2}\simeq1.301>1$

~

with eigenvalues $\lambda(A_{2})=\left\{ 1.2,\,1.3\right\} $ 

~

$\left\Vert A_{3}\right\Vert _{2}\simeq1.201>1$

~

with eigenvalues $\lambda(A_{3})=\left\{ 1.2,\,1.1\right\} $ 

So, the subsystems ruled by $A_{1}$ is stable, whereas the subsystems
ruled by $A_{2}$ and $A_{3}$ are unstable.

\begin{figure}
\includegraphics[width=0.55\textwidth]{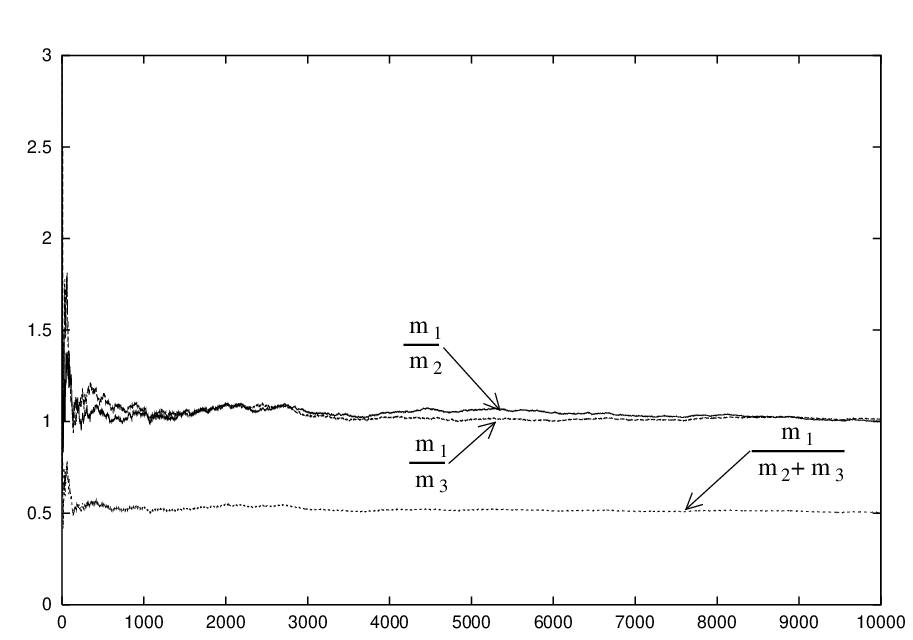}

\caption{\label{fig:3}Temporal evolution of $m_{1}/m_{2}$ , $m_{1}/m_{3}$
and $m_{1}/(m_{2}+m_{3})$.}
\end{figure}

After 10000 switches we obtain (see Figure \ref{fig:3})

~

\begin{equation}
\frac{m_{1}(t)}{m_{2}(t)}\underset{t\rightarrow\infty}{\rightarrow}k_{1,2}=1.0158>1\label{eq:ex_2.1}
\end{equation}

\begin{equation}
\frac{m_{1}(t)}{m_{3}(t)}\underset{t\rightarrow\infty}{\rightarrow}k_{1,3}=1.0236>1\label{eq:ex_2.2}
\end{equation}

with

~

$\left\Vert A_{1}\right\Vert _{2}\left\Vert A_{2}\right\Vert _{2}\left\Vert A_{3}\right\Vert _{2}\simeq0.973\leq1-\varepsilon_{\alpha}\quad\varepsilon_{\alpha}>0$

\begin{figure}

\begin{tabular}{cc}
{\footnotesize (a)}\includegraphics[width=0.55\textwidth]{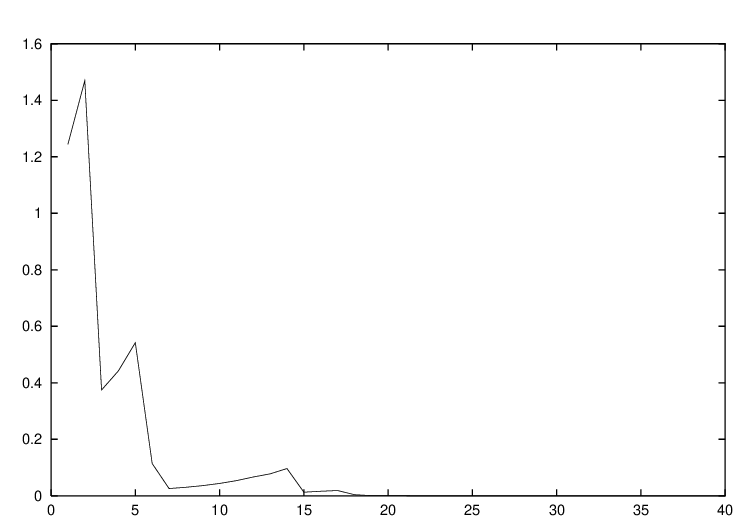} & {\footnotesize (b)}\includegraphics[width=0.55\textwidth]{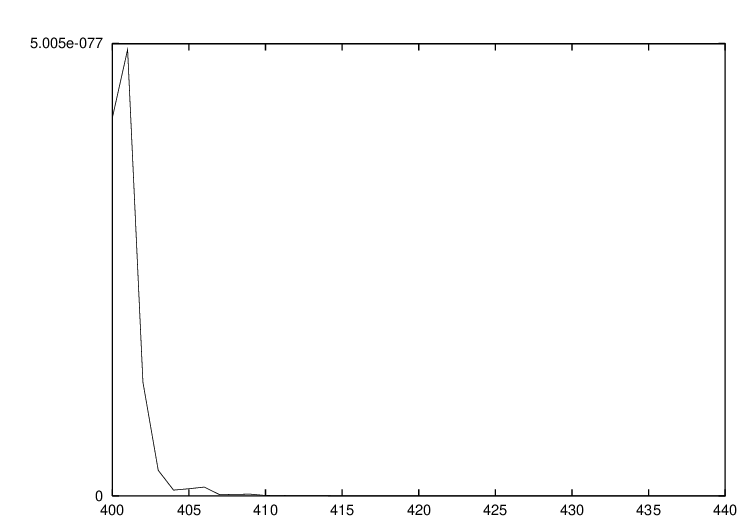}\tabularnewline
\end{tabular}

\caption{\label{fig:4}Temporal evolution of the asymptotic convergence in
norm to zero of discrete system (\ref{eq:ex_2}).}
\end{figure}

Given that the conditions of theorem 2 are satisfied the system must asymptotically
tend to zero, as it is shown in Figures \ref{fig:4}a and \ref{fig:4}b.

From (\ref{eq:ex_2.1}) and (\ref{eq:ex_2.2}) it results that

\[
\frac{m_{2}(t)+m_{3}(t)}{m_{1}(t)}\underset{t\rightarrow\infty}{\rightarrow\sim2}
\]
(as it can be seen in Figure \ref{fig:3}) that is, the system jumps
twice as many times into unstable subsystems as it does into stable
subsystems. As the switching sequence is generated by random numbers
it results, by using the Law of large numbers, that the system spends
twice as much time in unstable subsystem as in the stable subsystem.
However the system tends asymptotically to zero!

\section{Hybrid system}

When the system has both continuous and discrete subsystems we have
a hybrid system.

A linear hybrid system can be described by equations
\begin{equation}
\overset{.}{x}(t)=A_{\sigma_{1}(t)}x(t)\quad\qquad\sigma_{1}(t):\left[0,\infty\right)\rightarrow\left\{ 1,\cdots,n\right\} \label{eq:3-A}
\end{equation}

\begin{equation}
x(t+)=B_{\sigma_{2}(t)}x(t-)\quad\qquad\sigma_{2}(t):\left[0,\infty\right)\rightarrow\left\{ 0, 1,\cdots,m\right\} \label{eq:3-B}
\end{equation}
Here $x(t\pm)$ denote the one-sided limits of $x$ at $t$. The matrix
$B_0$ is the identity, and the signal $\sigma_2(t)$ that controls
discrete jumps has $\sigma_2(t)=0$ except at a discrete set of times
$t$.  The system evolves according to the differential equation
\eqref{eq:3-A}, except when the switching signal $\sigma_1(t)$
jumps (also a discrete set of times) or when $\sigma_2(t)>0$. 

These systems are more and more frequent in industry due to integration
of continuous and discrete systems. The continuous system might have
its origin in the flow or process of a factory or traffic, and the
discrete one in the digital control of the diverse steps of the process.
Hybrid systems show the same problems formulated by Liberzon y Morse,
that we have already mentioned formerly \citep{Sun2005a}. We can
deduce a theorem for these systems that gives the sufficient conditions
of stability by using theorems 1  and 2.

\paragraph{Theorem 3 \label{thm:3}}

Suppose a hybrid system given by (\ref{eq:3-A}) and (\ref{eq:3-B})
is such that the continuous subsystem (\ref{eq:3-A}) satisfies the conditions
of theorem 1  and the discrete subsystem (\ref{eq:3-B}) satisfies
the conditions of theorem 2. Then the hybrid system is asymptotically
stable.

\paragraph{Proof}

It is straight forward. It is enough taking norms of state $x(t)$
after a time $t$, and then to gather separately the norms corresponding
to the continuous subsystem and the discrete one. Then proofs of theorem
1  and 2  are respectively repeated for any group.$\blacksquare$

\paragraph{Remark 2\label{rem:thm3}}

If one of the subsystems has a bounded solution and the another one
tends asymptotically to zero (because it satisfies its respective
theorem) then the solution of hybrid system also tends asymptotically
to zero. We will return to this remark later.

\paragraph{Example 3}

\begin{figure}
\begin{tabular}{cc}
{\footnotesize (a)}\includegraphics[width=0.55\textwidth]{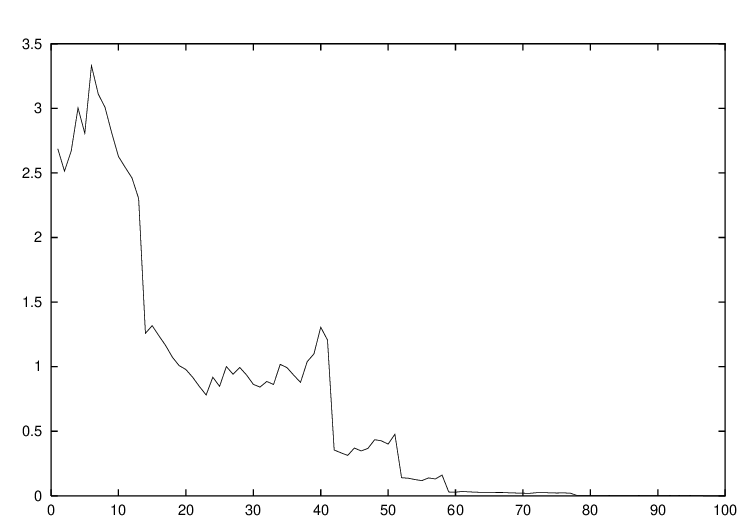} & {\footnotesize (b)}\includegraphics[width=0.55\textwidth]{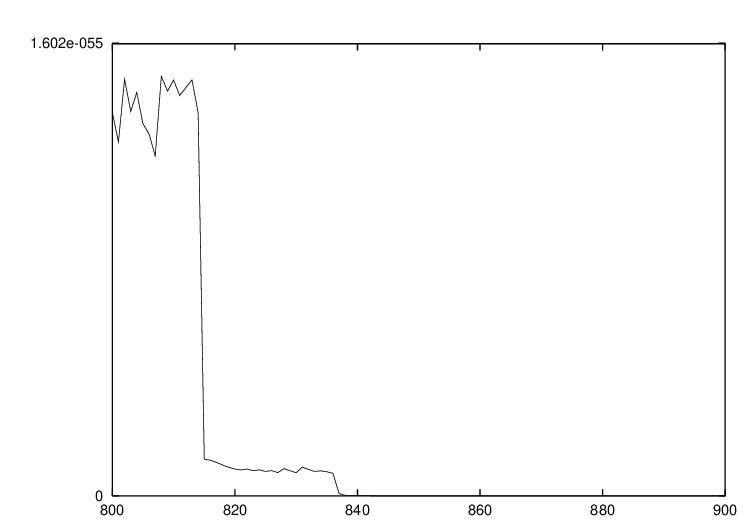}\tabularnewline
\end{tabular}

\caption{\label{fig:5}Temporal evolution of the asymptotic convergence in
norm to zero of the hybrid system (\ref{eq:3-A})-(\ref{eq:3-B}).}
\end{figure}

Finally, we consider the hybrid system \ref{eq:3-A} and \ref{eq:3-B}
formed by the continous one given in example 1 and the discrete system
given in example 2. Furthermore the activation of continuous and discrete
subsystems is given by a switching sequence $\sigma_{3}(t)$ generated
by random numbers. As theorem conditions are satisfied the system
tends asymptotically to zero as it is shown in Figures \ref{fig:5}a
and \ref{fig:5}b.

\section{Aftermaths: stabilisation via the control of the switching signal}

If the switching signal of a switched linear systems is not fixed, but
depends on a parameter or can be designed by the engineers, then
theorems in this paper allow the design of appropriate feedback
laws to make the system stable. Let us show how to do that.
\begin{description}
\item [{Continuous~time.-}] The condition $(i)$ of theorem 1 shows every
subsystem, described by $A_{i}$, must be left before a time $T_{i}$
(see remark 3  below).
\end{description}
Given that
\[
\left\Vert e^{A_{i}t}\right\Vert \leq p_{i}(t)e^{\mu_{i}t}
\]
where $\mu_{i}=\max\left\{ \mbox{Re}\lambda_{i}/\lambda_{i}\quad\mbox{eigenvalue}\quad\mbox{of}\quad A_{i}\right\} $
and $p_{i}(t)$ a polynomial of degree the order of $A_{i}$. If we
bound $p_{i}(t)\leq k_{i}$ in $\left[0,T_{i}\right]$ it follows
that
\[
\left\langle e^{A_{i}}\right\rangle \equiv\left[\prod_{i_{j}=1}^{m_{i}}\left\Vert e^{A_{i}t_{i_{j}}}\right\Vert \right]^{^{\frac{1}{m_{i}}}}\leq k_{i}e^{\mu_{i}\overline{t_{i}}}
\]
where
\[
\overline{t_{i}}=\frac{\sum_{i_{j}=1}^{m_{i}}t_{i_{j}}}{m_{i}}
\]
 is the average time that system (\ref{eq:(1)}) stays in subsystem
given by $A_{i}$.

Therefore
\begin{equation}
\prod_{i=1}^{n}\left\langle e^{A_{i}}\right\rangle \leq ke^{\sum_{i=1}^{n}\mu_{i}\overline{t_{i}}}\label{eq:F1}
\end{equation}

Thus, the time $\overline{t_{i}}$ can be deduced such that theorem
1  is satisfied and asymptotic stability is obtained. It is plain
to see that $\sigma(t)$ will not be unique, because we have an average
time $\overline{t_{i}}$, that is, engineers can choose any $\sigma(t)$
on the assumption that the average time $\overline{t_{i}}$ satisfies
theorem 1 .

It does not matter whether the switched linear system has unstable
matrices, the engineer must design the system in such a way that
it spends enough time (according to (\ref{eq:F1})) in stable matrix so
that they control the unstable matrices.

\paragraph{Remark 3\label{rem:9}}

The condition of staying a maximum $T_{i}$ in subsystem $A_{i}$
makes a lot of sense from a practical point of view. If the system
could be indefinitely in any subsystem then two possibilities would
arise:
\begin{enumerate}
\item Either the system is ruled by an unstable subsystem and it would be
destroyed.
\item Or the system is ruled by a stable subsystem and it would tend
asymptotically to zero. Then the result is trivial.
\end{enumerate}
So, the condition is not really a restriction at all.
\begin{description}
\item [{Discrete~time.-}] In this case, roughly speaking, theorem 2 
shows that the system must spend more time evolving under the stable
subsystem matrices than under the unstable ones.
\item [{Hybrid~system.-}] The system may have a stable solution because
both the continuous and discrete subsystems converge to 0. Or because
the continuous (discrete) subsystem converges to 0 while the discrete
(continuous) has a bounded solution (see remark 2).
\end{description}

\section{Conclusions}

We have presented theorems showing that, if a switched linear system
has a switching signal such that:
\begin{enumerate}
\item [i)]The geometric average of the subsystems flow norms is finite;
\item [ii)]The geometric averages of stable subsystems dominate the unstable
ones;
\item [iii)]The stable subsystems control unstable subsystems pair by pair;
\end{enumerate}
then the solution of the system converges asymptotically to zero (even
if the system spends more time in unstable subsystems ---see example 2).
The conditions are for continuous, discrete-time or hybrid systems.
These results would allow a practitioner to design the switching signal
in order to stabilize the system.

The engineer finds it easier to count how many times the system jumps
between different subsystems (see iii) than to calculate an average dwell
time switching in order to control the dangerous unstable subsystems.
It is easier, and what is more important in practical systems: it
is cheaper. Several numerical examples have been presented in order
to demonstrate the effectiveness of the theoretical findings. 


\end{document}